\documentclass{article}
\usepackage{amsfonts,amssymb,latexsym,xspace}
\usepackage{amsmath,enumerate}
\numberwithin{equation}{section}

\usepackage{pifont}

\newcommand\ItOsc{{\omega}}

\newcommand\nuh{\hat{\nu}}
\newcommand\Rh{R}

\newcommand\Q{{\mathcal Q}}
\newcommand\eps{\varepsilon}

\newcommand\muh{{\hat{\mu}}}

\renewcommand\L{{\mathcal L}}

\renewcommand\P{\mathcal P}

\newcommand\LL{{{\mathbf L}(d_0)}}

\newcommand\diam{\hbox{diam}}

\newcommand\om{\omega}

\newcommand\cte{\hbox{const}}
\newcommand{\un}{{\mathbf 1}}

\def\bin_#1^#2{\left(\smallmatrix #1 \\ #2 \endsmallmatrix\right)}

\newcommand{\RR}{\ensuremath{\mathcal{R}}}

\newcommand{\NN}{\ensuremath{\mathbb{N}}}
\newcommand{\N}{\ensuremath{\mathbb{N}}}

\newcommand\BP{{\mathcal B}}

\renewcommand{\H }{\"}

\newtheorem{theo}{\sc Theorem}[section]
\newtheorem{prop}[theo]{\sc Proposition}
\newtheorem{coro}[theo]{\sc Corollary}
\newtheorem{lem}[theo]{\sc Lemma}

\newtheorem{remnonum}{Remark}

\newenvironment{equa}
{\refstepcounter{equa} \ifhmode\par\medskip\fi \ \hfill
\begin{math} \displaystyle {\end{math}    \hfill (\theequa)
\medskip \\ }}

\newcounter{equa}
\renewcommand{\theequa}{\fnsymbol{equa}}
\newenvironment{dem}
{\ifhmode\par\fi \rmfamily \noindent {\it Proof\/}: }{\hfill
$\blacksquare$  \par\medskip \noindent}

\newenvironment{sketchdem}
{\ifhmode\par\fi\rmfamily \noindent {\it Sketch of proof \/}:
}{\hfill $\blacksquare$ \par \noindent}

\def\begindem{\par \noindent {\it Proof}\/:}
\def\eop{\hfill \ \qquad \qquad \blacksquare}
\title{Decay of correlations on towers with non-H{\"o}lder Jacobian and
  non-exponential return time}
\author{J{\'e}r{\^o}me Buzzi${}^1$ \\
V{\'e}ronique Maume-Deschamps${}^2$ \\
\\
${}^1$ Centre de Math{\'e}matiques de l'Ecole Polytechnique \\
${}^2$ Universit{\'e} de Bourgogne / Laboratoire de Topologie \ }
\date{April 2003}

\begin{document}
\maketitle
\footnote{\noindent {\em AMS classifiction 1991:} Primary 58F11 \\
{\em key words:} absolutely continuous invariant measures; equilibrium states; decay of correlations; transfer operator;  tower extension.}                      

\begin{abstract}
\noindent We establish upper bounds on the rate of decay of
correlations of tower systems with summable variation of  the
Jacobian and integrable return time. That is, we consider
situations in which the Jacobian is not H{\"o}lder and the return
time is only subexponentially decaying. We obtain a subexponential
bound on the correlations, which is essentially the slowest of the
decays of the variation of the Jacobian and of the return time.
\end{abstract}

\section*{ Introduction }

\noindent In this paper we study the speed of mixing, more precisely the rate of decay of correlations, of tower systems, a special class of countable Markov systems which naturally arise in the study of many dynamical systems by the procedure of {\bf induction} -- see \cite{Y1}.  Our goal is to provide a comprehensive statement in the following sense. There are two sources of loss of exponential speed: large return times and bad smoothness. By extending cone techniques, we deal simultaneously with both difficulties whereas previous works on decay of correlations \cite{KMS,BFG,Po,Y1} considered only one of these two obstructions. We prove that, although the analysis becomes more difficult when both obstructions are present, they nevertheless operate independently: the speed is just the minimum of the speeds allowed 1) by the defect in smoothness if the statistics of return times were exponential; 2) by the statistics of return times if we had H\"older smoothness.\\
\\
Let us state informally a corollary of our result:
\begin{theo}
Consider a tower system $F$ with a mixing invariant probability
measure $\muh$. Assume that the oscillation of the Jacobian on
$n$-cylinders is bounded by $n^{-\alpha}$ and the probability of
return time $n$ decays like $n^{-\beta}$.  Then, for sufficiently
smooth observables, the rate of decay of correlations is:
$$
     \left| \int \phi\cdot \psi\circ F^n \, d\muh -
            \int \phi\, d\muh \int \psi \, d\muh\right|
       \leq C \cdot K(\phi) \|\psi\|_{L^1} \cdot
                      \frac{1}{n^{\min(\alpha,\beta-\eps)-1}},
$$
for any $\eps>0$. $K(\phi)$ is some finite number depending only
on $\phi$~; $\|\psi\|_{L^1}$ is the $L^1$ norm w.r.t. the
reference measure.
\end{theo}

\noindent{\it Remarks.} \\
1. Our result allows returns which are not onto, which is quite convenient for applications.\\
2. The fact that the above bound depends on $\psi$ only through its $L^1$-norm is important for the study of asymptotic laws of return times \cite{C,CGS,Pa}.\\
3. The loss of $\eps$ in the exponent is probably due to our
method (in the H\"older-continuous case ($\alpha=\infty$ so to
speak), L.-S. Young \cite{Y1} obtained $\mathcal
O(n^{-\beta-1})$).

\medbreak

\noindent{\it An application to non-H\"older maps with an indifferent fixed point.} \\
Fix $0<\gamma<1/2$ and $\alpha>1$ and consider the interval map
$f:[0,1]\to[0,1]$ defined by $f(x)=
2^{1+\gamma}(x+x^{1+\gamma})/(2^\gamma+1)$ for $x<1/2$ and $f(x)=
\frac32(x-\frac12)+\frac{(x-\frac12)(\log2)^\alpha}
{2|\log(x-\frac12)|^\alpha}$ for $x>1/2$. Our result implies a
rate of correlation in $\mathcal
O(n^{-\min(\alpha,1/\gamma-\eps)+1})$ for arbitrarily small
$\eps>0$. \footnote{Indeed, using the notations of Section
\ref{sec:res}, one considers the tower defined by
$\Delta_0=[\frac12,1]$, the map $f_0(x)=f^{R(x)}(x)$ with the
return function $R(x)=\min\{n\geq1:f^n(x)\geq\frac12$ and
$\#\{k<n: f^k(x)\geq\frac12\}\geq\eps_0 n\}$ for some small
$\eps_0>0$. Then one can prove that $\omega_n=\mathcal
O(n^{-\alpha})$ and $\nu(\Delta_n)= \mathcal O(n^{-1/\gamma})$ and
apply our main theorem.} Our approach is the first to our
knowledge to be able to treat such maps.

\medbreak
\noindent Section \ref{sec:res} contains the  precise statement of
our results. We briefly recall  definitions and properties  of
Birkhoff's cones and projective metrics (section \ref{sec:cones})
and the construction of the a.c.i.m., establishing regularity of
the invariant density (section \ref{sec:mixing}). We define a
sequence of cones $C_j$ of ``Lipschitz'' functions (w.r.t. to an
ad-hoc metric) in section \ref{sec:defcones} and then establish
that the transfer operator iterated some $k_j$ times sends one
cone into the next by a $\gamma_j$-contraction in  section
\ref{sec:contraction} for some semi-explicit $\gamma_j<1$. Finally
in section \ref{sec:conclusion}, we deduce from this a convergence
in the uniform norm at speed $\prod_{p=1}^{j}\gamma_p$, with $j$
largest such that $k_1+\dots+k_j\leq n$, and make this estimate
explicit in the exponential, stretched exponential and polynomial
cases. \medbreak

\noindent {\bf Acknowledgments:} The authors are grateful to {\it
program ESF/PRODYN} which has partially supported  the {\it
International Conference on Dynamical Systems, Abbey of ``La
Bussi{\`e}re''} where part of this work was carried through.
\section{Setting, statement of the results}\label{sec:res}
Let us describe our tower model which follows Young's \cite{Y1}.
A tower is defined by:
\begin{itemize}
\item a basis, which is a probability space $(\Delta_0,m_0)$
together with a non-singular self-map $f_0$; \item a partition
$\Delta_{0,j}$, $j\in\N$ such that    $f_0:\Delta_{0,j}\to
f_0(\Delta_{0,j})$ is one-to-one and satisfies
$f_0(\Delta_{0,j})$ is a union of some $\Delta_{0,k}$, for some
$k$'s; \item a return time, i.e., a function $R:\Delta_0\to\N$,
constant on each $\Delta_{0,j}$, $j\in\N$.
\end{itemize}
The tower $\Delta$ is then the disjoint union of the floors
$\Delta_\ell$, $\ell \in \N$:
$$
 \Delta_\ell = \{ (x, \ell) \ | \ x \in \Delta_0 \/, \ R(x) > \ell \/ \} \/.
$$
It is endowed with the measure $\hat\nu$ which is just the
restriction of the copy of $m_0$ on each floor.   We will denote
by $\Delta_{\ell\/,j}$, $ \ell < R_{|\Delta_0^j}$ the copy of
$\Delta_{0\/,j}$ inside $\Delta_\ell$:
$$
\Delta_{\ell\/,j} = \{ (x, \ell) \ |\  x \in \Delta_{0\/,j} \/, \
R(x) > \ell \/ \}   \/.
$$

\noindent The dynamic on the tower, $F~: \Delta \longrightarrow
\Delta$, is defined  by
$$
\left\{
\begin{array}{lcl}
F(x,\ell) & = & (x, \ell +1) \ \mbox{if } R(x) > \ell +1 \\
      & = & (f_0(x),0)  \ \mbox{otherwise} \/.
\end{array}
\right.
$$
One can think of $F$ as the unfolding of the underlying induction:
in applications, $F$ will be often conjugate to the original map,
which $f_0$ is some (variable) power.  We assume that the
partition $\RR = \{ \Delta_{\ell\/,j} \}$ generates in the sense
that $\displaystyle \bigvee_{i=0}^{\infty} F^{-i} \RR $ is the
partition into points mod $\hat\nu$. For $k \in \N$, the elements
of the partition $\displaystyle \RR^{(k)} = \bigvee_{i=0}^{k-1}
F^{-i} \RR$ are called {\bf cylinders or   $k$-cylinders}. We
denote by $C_k(x)$ the element of $\RR^{(k)}$ which contains $x$.
Let $JF$ be the Jacobian of $F$ with respect to $\nuh$ (this
Jacobian is well defined because of the non singularity of $f_0$).
The modulus of continuity of $JF$ will be controlled by the
following, dynamics-dependent sequence:
$$
    \ItOsc_n  = \sup_{C\in\RR^{n}} \sup_{x,y\in C} \log\frac{JF(x)}{JF(y)}.
$$
For $x, y \in \Delta$ the {\bf separation time} $s(x,y)$ is the
largest integer $n \geq 0$ such that for all $0\leq j \leq n$,
$F^j (x)$ and $F^j (y)$ belong to the same atom of the partition
$\RR$. Set
$$
     d_0 (x,y) =\sum_{j \geq s(x,y) +1} \ItOsc_j   \/.
$$
Note that the metric $d_0$  is designed so that the family of
functions:
$$
\log JF^n = \log \prod_{i=0}^{n-1} JF \circ F^i.
$$
are uniformly Lipschitz w.r.t.\ it.\\
Let us summarize our assumptions on  the tower.
\begin{itemize}
\item[(A.I)] {\bf Summability of upper floors.}
$$
 \sum_{\ell \in \N} \hat{\nu} \left(\{ x \in \Delta_0 \ | \ R(x) > \ell \}\right) = 1\/.
$$
\item[(A.II)] {\bf Generating Partition.} The partition $\RR$
generates under $F$ i.e.: the partition $\bigvee_{n=0}^\infty
F^{-n} \RR$ is the partition into points. In particular, $d_0$
defines a metric on $\Delta$.
\item[(A.III)] {\bf Summable variation.} Let $JF$ be the Jacobian
of
  $F$ with respect to $\nuh$.  We assume that $JF$ satisfies:
$$
\sum_{n \in \NN} \ItOsc_n  < \infty.
$$
\item[(A.IV)] {\bf Large image and Markov properties.} Each
$F^{\Rh}\Delta_{0,j}$ is a union of some $\Delta_{0,p}$,
$p\in\NN$, ({\bf Markov property}) and  ({\bf Large image}):
$$
    \eta := \inf_{j \in \NN} \hat{\nu} (F^{\Rh}(\Delta_{0,j})) >0.
$$
\end{itemize}
Contrarily to \cite{Y1} we do not assume the Bernoulli property:
$f_0(\Delta_{0,j})=\Delta_0$, but only the weaker Markov property
above. The collection of sets $f_0(\Delta_{0,j})$ defines a
partition $\mathcal B$ which is less refined than
$\{\Delta_{0,j}\}_{j\in\N}$, so that it is in particular countable
$\mathcal B=\{B_1,B_2,\dots\}$. Remark that, by an easy induction,
if $x,y$ are contained in the same element of $\mathcal B$, then
the pre-images of all orders of $x$ and $y$ are paired in the
following sense. \medbreak

Given $x,y\in\Delta_0$, say that $x',y'\in\Delta$ are {\bf paired
pre-images} if $F^n x'=x$, $F^n y'=y$ and $F^k(x')$ nd $F^k(y')$
belong to the same element of $\RR$ for all $0\leq k<n$. Observe
that (A.III) implies that in this situation we have:
\begin{equation}
     \left\vert \frac{JF^n (x')}{JF^n (y')} - 1\right\vert \leq  C \cdot d_0(x,y) \/, \ \qquad{\rm with\ \ } C=\exp\sum_{j\geq1} \ItOsc_j \/. \label{bd}
\end{equation}
This is ``bounded distortion''.
\begin{remnonum}
In \cite{BM}, we proved that multi-dimensional piecewise expanding
maps in higher dimension are (under quite general hypothesis)
conjugate to such a tower map.
\end{remnonum}
Let $\LL$ be the space of bounded functions on $\Delta$ that are
{\bf locally Lipschitz} with respect to the metric $d_0$, i.e.,
for some $K<\infty$, for all  $x$, $y$ in the same $B_{j,\ell}$,
$$
 \vert \varphi (x) - \varphi (y)\vert \leq K d_0(x,y).
$$
$K(\varphi)$ is the smallest number $K$ such that the above inequality is satisfied. Let $\Vert \varphi \Vert_\LL = K(\varphi) + \ \Vert \varphi \Vert_\infty$ be the norm on $\LL$.\\
\\
To study the ergodic properties of $F$, we have to decompose it
into topologically mixing components. Observe that  $\RR$  has a
natural graph structure: $P\to Q$ iff $F(P)\supset Q$. Its
(restricted) {\bf spectral decomposition}  is ${\cal P}={\cal P}_t
\cup \bigcup_{i} \bigcup_{j=0}^{p_i-1} {\cal R}^{(i)}_j$,  where:
\begin{itemize}
 \item ${\cal P}_t$ is the set of {\bf transient} elements of $\cal P$, i.e., elements $P$ such that there exists a  path from $P$ going to some $Q\in\cal P$ and there is no path from $Q$ to $P$ (observe that we don't decompose this part into irreducible subsets). The elements that are not transient are called {\bf recurrent}. \item for each $i$, $\bigcup_{j=0}^{p_i-1} {\RR}^{(i)}_j$ is the set of $P\in\RR$ such that there exist paths from $P$ to $Q$ and $Q$ to $P$, for some fixed $Q=Q(i)$ (i.e., these unions are the irreducible components of $\RR$  from which no arrows leave).
 \item if there is an arrow from ${\RR}^{(i)}_j$ to ${\RR}^{(k)}_l$ then $k=i$ and $l=j+1 \ \mbox{mod} \ p_i$.
\end{itemize}
Finally, $\Delta^{(i)}_j$ is the union of the elements of ${\RR}^{(i)}_j$. Observe that, up to  trivialities, it is enough to study the dynamics of $F^{p_i}:\Delta^{(i)}_0\to \Delta^{(i)}_0$ for each $i$. We call this  the {\bf spectral reduction}.\\
\\
Our main result is the following theorem.
\begin{theo}\label{theotour}
Let $(\Delta, F, \nuh)$ be a tower system satisfying (A.I - IV).
First, there exists an  invariant probability measure absolutely
continuous with respect to $\nuh$ (a $\nuh$-a.c.i.m. for short).
\\
Second, any $\nuh$-a.c.i.m.\ $\muh$, up to the spectral reduction,
is mixing, with the following speed estimate: for all
$\varphi\in\LL$ and $\psi\in L^\infty(\Delta)$,
$$
   \left| \int\limits_\Delta \varphi\circ F^n \cdot \psi \, d\muh -\int\limits_\Delta \varphi \, d\muh\int\limits_\Delta \psi \, d\muh \right|   \leq C \cdot \Vert \psi\Vert_\LL \/  \|\varphi\|_{L^1(\muh)} \cdot u_n \qquad \hbox{for all }n\geq0
$$
for some $C<\infty$ and a sequence $u=(u_n)_{n=0}^\infty$
converging to zero which can be made explicit:
\begin{itemize}
\item if $\ItOsc_n = O (\rho^n)$ for some $0<\rho<1$ and $\nuh
(\Delta_n) =O(\alpha^n)$  for some $0< \alpha < 1$  then $u_n
=\kappa^n$ for some $0<\kappa<1$, \item if $\ItOsc_n = O
(n^{-\alpha})$ for some $\alpha >1$ and   $\nuh (\Delta_n)
=O(n^{-\beta})$ for some $\beta >1$ then  $u_n =n^{-\min
\left(\alpha, \beta-\eps\right)-1}$ for all $\eps>0$. \item if
$\ItOsc_n = O (e^{-n^\alpha})$ and $\nuh (\Delta_n)
=O(e^{-n^\beta})$ for some $0<\alpha,\beta<1$, then $u_n =
e^{-n^{\min(\alpha,\beta)-\eps}} $ for all $\eps>0$.
\end{itemize}
\end{theo}

\medbreak

\section{Birkhoff's cones and projective metrics}\label{sec:cones}
The main tool for the proof of  Theorem~\ref{theotour} will be the
theory of cones and projective metrics of Garrett~Birkhoff
\cite{Bi1}. P.  Ferrero and B.~Schmitt \cite{FS} applied it to
estimate the correlation  decay for random products of matrices.
Recently this strategy has been used by many authors to obtain
exponential decay of correlations (see for example \cite{L}).  We
are closer to \cite{KMS} and \cite{M} which have used these
techniques in a different way to  obtain sub-exponential decay of
correlations.  Let us recall definitions  and properties of cones
and projective metrics (see \cite{L}  for  a more complete
presentation). Let  $B$ be a  vector space and let   $C \subset B$
be a {\bf Birkhoff
  cone}, i.e., a cone  with the following  properties.
\begin{itemize}
\item $C$ is convex,
\item $C \/ \cap \/ -C = \{0\}$,
\item if $\alpha_n$ is a sequence of real numbers such that
$\alpha_n \rightarrow \alpha$ and $x - \alpha_n y \in C $ for all
$ n$, then $x -\alpha y \in C$. This property is called ``integral
closure''.
\end{itemize}
Such a cone is endowed with the pseudo-metric $\delta_C$ on $C$
defined in the following way (it is pseudo because it is not
necessarily finite and it does not separate points). For $x, y \in
C$,
$$
        \mu (x,y) = \inf \{ \beta>0 \ \mbox{such that} \ \beta x - y \in C \}.
$$
with the convention: $\mu(x,y) = \infty$
if the corresponding set is empty.  Let $\delta_{C}(x,y) = \log \mu(x,y)\mu(y,x)$. We remark that $\delta_C$ satisfies the triangle inequality: if $\beta x - y \in C$ and $\tilde{\beta} y - z \in C$ then $\beta \tilde{\beta} x - z \in C$ since $C$ is a convex cone, so $\mu(x\/, z) \leq \mu(x\/,y)\cdot \mu(y\/,z)$ and the triangle inequality follows.   Finally, observe that $\delta_C$  is projective: $\delta(x,y)=0 \iff $ $x$ and $y$ are colinear.  \\
\\
The usefulness of this projective metric is that it allows a
`geometric' proof of the contraction through the following result.
\begin{theo} {\rm \cite{Bi1}} \label{projm1}
Let $C$ and $C'$ be two cones, $P$ a linear operator  $P : C
\rightarrow C'$. Let $\Gamma$ denote the diameter of  $PC$ in
$C'$:
$$
 \Gamma = \sup_{f, g \in C} \delta_{C'}(Pf, Pg) \leq \infty.
$$
For any $f, g$ in $C$, we have:
$$
 \delta_{C'}(Pf, Pg) \leq \tanh \left(\frac{\Gamma}{4} \right) \/\delta_C(f,g) \/.
$$
\end{theo}

\noindent  This theorem implies that a linear map between cones never increases distances and is in fact a contraction as soon as $\Gamma < \infty$.  \\
The following result allows the translation of contraction w.r.t.
cone metric to a contraction w.r.t. norm. A norm  $\Vert \ \Vert$
on  $B$ is {\bf adapted to $C$}, if for $f$ and $g$ in  $B$ such
that both $f+g$ and $f-g$ belong to $C$, then $\Vert g \Vert \leq
\Vert f \Vert$.  $\rho:C\to\RR_+$ is a {\bf homogeneous form
adapted to   $C$} if, i) for any  $\lambda >0$~; ii) $f \in C$,
$\rho(\lambda f) = \lambda \rho(f)$  and if  $f-g \in C$ implies
$\rho(g) \leq \rho(f)$.
\begin{theo} {\rm \cite{Bi1}, \cite{L}.} \label{projm2}
Let $C$ be a Birkhoff cone,   let $\Vert \ \Vert$ and  $\rho$ be
adapted to  $C$. For any  $f$ and $g$ in $C$ such that $\rho (f) =
\rho(g) \neq 0$ we have:
$$
\Vert f- g \Vert \leq (e^{\delta(f,g)} -1) \min (\Vert f \Vert ,\/
\Vert g \Vert) \/.
$$
\end{theo}

\section{Construction of a $\nuh$-a.c.i.m.}\label{sec:mixing}

As usual,  the {\bf transfer operator} acting on bounded functions
is defined by:
$$
 \L_0 f (x) = \sum_{Fy = x} \frac{1}{JF(y)} f(y) \/.
$$
The measure $\nuh$ is conformal for $\L_0$ in the following sense:
for any bounded function $f$,
$$
 \int \L_0 f d\nuh = \int f d\nuh \/.
$$
For $s \in \NN$, the {\bf $s$-cylinders} are the non empty sets of the form: $\bigcap_{i=0}^{s-1} F^{-i} A_{i}$ with $A_{i} \in \RR$. For $k \in \NN$ and $x \in \Delta$, $C_k(x)$ denotes the $k$-cylinder which contains $x$.\\
The following lemmas are technical tools to study $\L_0^n$.\\
\\
\begin{lem}\label{BD}
There exists $C<\infty$ such that for any $\ell \in \NN$ and any
$x \in \Delta_\ell$ and $k \in \NN$ with $F^{k} x \in \Delta_0$,
$$
 C^{-1} \nuh (C_k (x)) \leq \frac{1}{JF^k (x)} \leq C \nuh (C_k (x)) \/.
$$
\end{lem}
\begindem
Let $x \in \Delta_\ell$ such that $F^k(x)\in\Delta_0$. The Markov
property and the large image property (A.IV) imply that $\nuh (F^k
C_{k} (x)) \geq \eta >0$. The bounded distortion property
(\ref{bd}) gives:
\begin{eqnarray*}
C^{-1} \frac{\nuh(C_k(x))}{\nuh(F^kC_k(x))} & \leq \displaystyle\frac{1}{JF^{k} (x)} \leq & C \frac{\nuh(C_k(x))}{\nuh(F^kC_k(x))} \\
C^{-1} \frac{\nuh(C_k(x))}{1} & \leq \displaystyle \frac{1}{JF^{k}
(x)} \leq & C \frac{\nuh(C_k(x))}{\eta}
\end{eqnarray*}
The Lemma is proved. $\eop$
\\
\noindent
\begin{lem}\label{bounded}
There exists $K<\infty$ such that:
\begin{itemize}
\item for all $x\in \Delta$, all $n \in \NN$, $ \L_0^n \un (x)
\leq K$. \item for all $x,y$ in a given $B_{j, \ell}$ and all $n
\in \NN$,:
\begin{equation}\label{distorsion}
\vert \L_0^n \un (x) - \L_0^n \un (y) \vert \leq K d_0(x,y) \/.
\end{equation}
\end{itemize}
\end{lem}
\begin{dem}
The upper  bound $\L_0^n \un \le K$ follows from Lemma \ref{BD},
by writing:
\begin{equation}\label{upper}
\L_0^n\un (x) = \sum_{x'\in F^{-n}x} \frac{1}{JF^n(x')} \leq C
\sum_{x'\in F^{-n}x} \nuh(C_n(x')) \leq C
\end{equation}
Let $x$ and $y$ belong to one $B_{j,\ell}$. Their preimages by
$F^n$ are paired, i.e., if $F^n x' = x$, there is exactly one $y'
\in C_n (x')$ such that $F^n y'=y$. So, using (A.IV), we get:
\begin{eqnarray*}
\vert \L_0^n \un (x) -\L_0^n \un (y) \vert & =& \left|  \sum_{F^n x' = x} JF^n (x')^{-1} - \sum_{F^n y' = y} JF^n (y')^{-1} \right|    \\
&  = &  \left|  \sum_{F^n x' = x} JF^n (x')^{-1}\left(\frac{JF^n(x')}{JF^n(y')} - 1\right) \right|     \\
 &  \leq & C \L_0^n\un(x) d_0(x,y)    \\
 &  \leq & KC d_0(x,y) \/
\end{eqnarray*}
(\ref{distorsion}) is proved.
\end{dem}
\begin{coro}
$F$ admits a $\nuh$ a.c.i.m.
\end{coro}
\begin{dem}
By Lemma \ref{bounded}, the sequence $\frac1n \sum_{i=0}^{n-1}
\L_0^i \un $ is relatively compact for the topology of uniform
convergence on compact subsets (this is Arzela-Ascoli theorem on
separable spaces). Each limit point $h$ of this sequence is a non
zero  fixed point for $\L_0$  (by Lebesgue's dominated theorem,
$\nuh (h) =1$), so that $\muh = h \nuh$ is a $\nuh$-a.c.i.m.
\end{dem}
The system $(\Delta, F, \nuh, \RR)$ has a Markov structure in the
sense that for each $P \in \RR$, $F(P)$ is a union of atom of
$\RR$. According to \cite{ADU}, we will say that $F$ is {\bf
aperiodic} if:
\begin{equation} \label{mixing}
\forall \ P \/, \ P' \in \RR \  \exists N \in \NN \ \mbox{such
that} \  \nuh (F^{-n} P \cap P') >0 \ \forall \ n \geq N \/.
\end{equation}
The existence of a $\nuh$-a.c.i.m. implies that the recurrent part
is non empty (it contains the support of $\muh$). Up to the
spectral reduction, {\bf we may and shall assume that $F$ is
aperiodic}. We remark that aperiodicity implies that  any
$s$-cylinder has positive $\nuh$-measure. The following lemma
implies that any $s$-cylinder also has $\muh$ positive measure.
\begin{lem}\label{inf}
If $F$ is aperiodic   then $h (x) >0$ for all $x \in \Delta$. Moreover, \\
$\inf \muh[F^k (C_k(x))] >0$ where the inf is taken on all $k \in
\N$ and $x$ such that $F^k x \in \Delta_0$.
\end{lem}
\begin{dem}
Theorems 2.5 and 3.2 in \cite{ADU} imply that if $F$ is aperiodic
then  $\L^n_0 \un \rightarrow h$ uniformly on compact sets. Let
$K$ be given by Lemma \ref{bounded}. We have for  $j = 1 \/,
\ldots \/, $ any $\ell \/, n \in \NN$, $x, y \in B_{\ell\/,j}$,
their paired preimages will be denoted by $x'$ and $y'$,
\begin{eqnarray}
 \L_0^n \un (x)& =& \sum_{F^n x' = x } JF^n(x')^{-1} =  \sum_{F^n y' = y } JF^n(y')^{-1} \/ \frac{JF^n(y')}{JF^n(x')}\nonumber \\
& \leq &(C d_0(x,y)+ 1) \sum_{F^n y' =y} JF^n(y')^{-1} \leq K
\L_0^n \un (y)  \ \mbox{using (\ref{bd})} \/. \label{supinf}
\end{eqnarray}
Taking the limit when $n$ goes to infinity implies: for $x, y \in
B_{\ell\/,j}$,
\begin{equation}\label{supinfh}
h(x) \leq K h(y)  \/.
\end{equation}
So, for all $(j\/, \ell)$, either $h \equiv 0$ on $B_{\ell \/, j}$
or $h >0$ on $B_{\ell \/, j}$. But $h|_{B_{j,\ell}}\equiv0$
implies that $\nuh(B_{j,\ell})=0$, a contradiction to the
aperiodicity. This concludes the proof of the first part of the
lemma.

To prove the second part, let us remark first that
$\nuh[F^kC_k(x)] \geq \eta >0$ for all $k$ and $x$ such that $F^k
x \in \Delta_0$. Also, the Markov property implies   that there
exists finitely many integers $i_1 , \/ \ldots \/, i_p$ such that
each $F^k C_k(x)$ contains at least one $\Delta_{0\/, i_j}$, $j=1
\/, \ldots \/, p$. This implies the announced result using that
$h>0$.
\end{dem}

\noindent Let us note that Lemma \ref{inf} implies that $\muh(P)
>0$ for any cylinder $P$. The following lemma  is a direct
consequence of mixing.
\begin{lem}\label{exactness}
There exists positive numbers $A$ and $B$ such that for any $f$,
$g \in L^2 (\nuh)$, with $\muh (f) >0 $, $\muh (g) >0$, there
exists $n_0$ such that for $n \ge n_0$,
$$
 A \le \frac{\muh (f \circ F^n \cdot g)}{\muh(f) \muh (g)} \le B \/.
$$
\end{lem}

We shall now construct a sequence of cones $C_j$ and a sequence of
integers $k_j$ such that $\L^{k_j}$ maps $C_{j-1}$ into $C_{j}$
with uniformly bounded diameter, where $\L$ is the normalized
transfer operator defined as follows:   $\L f = \frac1h \L_0
(fh)$. Because of Lemma \ref{inf}, $\L$ is well defined. Moreover
it satisfies: $\L \un = \un$. The Jacobian of $F^n$ with respect
to $\muh$ is:
$$
\frac{JF^n \cdot h \circ F^n}{h}
$$
Let $x$ and $y$ belong to the same $B_{ \ell \/,j}$, $x'$ and $y'$
be their paired preimages by $F^n$.  Following eqs. (\ref{supinf}
- \ref{supinfh}), we get, for some $C'>C$:
\begin{eqnarray*}
(1-C'd_0(x'\/, y') ) &\leq \frac{h(x')}{h(y')} \leq & (1+C'd_0(x'\/, y')) \/, \\
(1-C'd_0(x\/, y) ) &\leq \frac{h\circ F^n (x')}{h\circ F^n (y')}
\leq & (1+C'd_0(x\/, y)) \/.
\end{eqnarray*}
We deduce  that the Jacobian of $F^n$ with respect to $\muh$
satisfies a bounded distortion inequality like (\ref{bd}) with an
appropriate constant that we will continue to denote by $C$. From
now on, we abuse notations and  $JF$ will be the Jacobian of $F$
with respect to the invariant measure $\muh$.  We remark that the
proof of Lemma \ref{BD} and Lemma \ref{inf}  give for some $C>0$:
\begin{equation}\label{gibbsmu}
C^{-1} \muh (C_k (x)) \leq \frac{1}{JF^k (x)} \leq C \muh (C_k
(x)) \/.
\end{equation}

\section{The cones} \label{sec:defcones}

\subsection{Auxiliary definitions}
In what follows $\muh$ is a mixing a.c.i.m. on $\Delta$.\\
\\
We need first some auxiliary definitions. We set for convenience $D=5$. \\
Let $(v_n)_{n\in\NN}$ be such that:
$$
     \sum_{n\geq1} v_n \cdot \nuh(\Delta_n) < \infty \text{ and }    v_n\to\infty.
$$
We also assume that $v_n$, and for each $k\in\NN$, $v_n/v_{n+k}$
are non-decreasing functions of $n$.  We define $\muh_v:=v \cdot
\muh$ where we have introduced the function $v = \sum_{\ell\geq0}
v_\ell \cdot 1_{\Delta_\ell}$.   Let $R_0(p) =\sum_{k>p}
\omega_k(g)$. We pick an integer $s$ so large that
$$
      R_0(s) \leq 10^{-5}.
$$
Let $P_\infty=\bigcup_{i\geq t\atop \ell\geq0} \Delta_{\ell,i}
\cup \bigcup_{\ell\geq t \atop i\geq0} \Delta_{\ell,i}$ with the
parameter $t$ chosen so large that:
$$
   \frac{\muh_v(P_\infty)}{\eta} \leq 10^{-5} \/,
$$
where  $\eta=\inf_j \nuh(F^R (\Delta_{0\/, j})) >0$.\\
Let $\Q_1$ be the finite collection of $s$-cylinders covering
$\Delta\setminus P_\infty$. Let $\Q$ be the finite partition of
$\Delta$ defined as $\Q_1\cup\{P_\infty\}$.   Let $k_0$ be such
that for all $k\geq k_0$, for all $P,Q\in\Q$:
$$
\begin{aligned}
       &\frac78 \leq \frac{\muh(F^{-k}P\cap Q)}{\muh(P)\muh(Q)}        \leq \frac98 \\
      &\frac78 \leq \frac{\muh_v(F^{-k}P\cap Q)}{\muh(P)\muh_v(Q)}        \leq \frac98 \\
 \end{aligned}
$$
Such a $k_0$ exists as $(F,\muh)$ is mixing, $\Q$ is finite and
the function $v$ is in $L^1(\muh)$.
\subsection{Definition of the distances $d_j$}
We set $d_j(x,y)=R_j(s(x,y))$ where the functions $R_j(\cdot)$ are
defined inductively in the following way. Recall that $R_0(\cdot)$
and $k_0$ have been defined above.   Assuming that
$R_{j-1}(\cdot)$ is defined we set:
$$
   k_j = \min \{ k\geq k_0 : R_{j-1}(s+k) \leq D^{-1} R_0(s) \}
$$
 and
 $$
  R_j(p) = D[ R_0(p)+R_{j-1}(p+k_j) ].
$$
We observe that, $\Q_1$ being a collection of $s$-cylinders, its $d_j$-diameter is bounded by $R_j(s)=D[ R_0(p)+R_{j-1}(p+k_\ell) ] \leq (D+1) R_0(s)$, a number independent of $j$.  \\
We introduce the auxiliary values $q(j)=k_1+\dots+k_j$.
\subsection{Definition of the cones}
As stated in the introduction, we are going to prove Theorem \ref{theotour} by cone techniques. Let us explain a bit how to construct the cones and how sub exponential decay of correlations may be obtained.  \\
\\
We start by recalling the classical way of using cones (see \cite
{FS} and \cite{L} for details). To get exponential decay of
correlations, it is sufficient to find a cone $C$ and an integer
$k$ such that $\L^k$ maps $C$ into itself and the diameter
$\Gamma$  of $\L^k C $ into $C$ is finite. If the fixed point $h$
of $\L$ belongs to $C$ then Theorem \ref{projm1} gives, for any
integer $j$:
$$
  \delta_C (\L^{kj} f \/, \ h) \leq \gamma^{j-1} \Gamma \ \mbox{where} \ \gamma = \tanh\frac{\Gamma}{4} < 1 \/.
$$
Hence Theorem \ref{projm2} gives that for $f \in C$, $\Vert \L^p f \ - \ h m(f)\Vert $ goes to zero exponentially fast for $\Vert \ \Vert$ an adapted norm, provided $f \mapsto m(f)$ is adapted. Then one has to extend this result from the cone to a Banach space.\\
\\
The starting point of the construction of cones is usually a Lasota-Yorke inequality (it will be done in section \ref{contrac2}). If the metric $d_0$ is not of exponential type (i.e., $d_0(x,y) \leq \beta^{s(x,y)}$ with $0<\beta<1$), then we cannot obtain a  Lasota-Yorke inequality. This is why we have to introduce the sequence of metric $d_j$ and a sequence of cones.  Roughly speaking, for any integer $j$, we will consider a cone $C_j$ of functions $f$ that are locally Lipschitz for the metric $d_j$ and the Lipschitz constant of which is controlled ({\bf see condition   2. below}). Thanks to the definition of the metric $d_j$ and to the Lasota-Yorke inequality, if $f$ belongs to $C_j$ then $\L^{k_j} f$ will be locally Lipschitz with respect to the metric $d_{j+1}$ and we will control its Lipschitz constant. This will imply that $\L^{k_j}$ maps $C_j$ into $C_{j+1}$ with finite diameter $\Gamma$. Then, using Theorems \ref{projm1} and \ref{projm2}, we get that for $f$ in $C_0$, $\L^{k_1 + \cdots + k_j}f$ goes to $m(f) h$ at rate $\gamma^j$ (with $\gamma = \frac{\tanh \Gamma}{4}$) in any adapted norm, provided that $h$ belongs to all the cones $C_j$ and that $f \mapsto m(f)$ is adapted. This is the philosophy of the construction.   \\
A source of difficulty is the following. To ensure that a cone $C$
satisfies properties of section \ref{sec:cones} and more
specifically the condition $C \cap - C  = \{0\}$,  some positivity
for the functions in the cone is needed. On the other hand, if $C
\subset \{ f \geq 0 \} =: C_+$ then for any $f$, $g$ in $C$,
$\theta_C (f,g) \geq \theta_{C_+} (f,g)$ (use Theorem \ref{projm1}
with $P=Id$) and
$$
  \theta_{C_+} (f,g) = \frac{\sup f}{\inf f} \cdot \frac{\sup g}{\inf    g} \/.
$$
Since the functions of the cones are only {\it locally} Lipschitz,
we will have a good control on $\frac{\sup f}{\inf f}$ on each
floor $\Delta_\ell$ but not on the whole space $\Delta$. Observe
that because of the definition of $\L$, we cannot hope to control
globally Lipschitz constant (just try to compute $| \L f(x) - \L
f(y) |$ for $x \in \Delta_0$ and $y \in \Delta_\ell$, $\ell >0$)
and we have to restrict ourselves to locally Lipschitz functions.
This problem is solved by considering the finite  partition $\Q$
of $\Delta$ which is decomposed into finitely many $s$-cylinders
(the ``compact'' part) and the complementary of the union of these
$s$-cylinders (the ``non compact'' part). Then we require the
positivity of some kind of conditional expectation of $f$ with
respect to this finite partition ({\bf see condition 1. below}).
This together with the control of the local Lipschitz constant
leads to a good control of $f$ on the atoms of the compact part.
Then, we requir!
 e !
!
another kind of control on the non compact part ({\bf see conditions 3   and 4  below}).  \\
\\
The cone $C_j(a,b,c)$ is the set of all real functions $f$ on
$\Delta$ satisfying the following conditions:
\begin{enumerate}
\item $a \cdot E_\muh(f) \leq E_\muh(f|\Q) \leq 6b\cdot
E_\muh(f)$. \item for all $x,y\in\Delta$ with $\BP(x)=\BP(y)$,
$$
    |f(x)-f(y)|\leq 12b\cdot E_\muh(f) \cdot d_j(x,y).
$$
\item for all $\ell\leq q(j)$,
$$
   \sup_{\P_\infty\cap\Delta_\ell} |f|            \leq 90c\cdot v_\ell \cdot E_\muh(f).
$$
\item for all $\ell>q(j)$,
$$
    \sup_{\P_\infty\cap\Delta_\ell} |f |  \leq 90 c \cdot v_{q(j)} \cdot E_\muh(f).
$$
\end{enumerate}

\section{Contraction of the cones}\label{sec:contraction}
The purpose of this section is to prove the following proposition.
\begin{prop}\label{image}
We have
 $$
   \L^{k_j} C_j(0,1,1) \subset  C_{j+1}\left(\frac45,\frac15,   \max\left(\frac15,\frac{v_{q(j)}}{v_{q(j+1)}} \right)  \right)
$$
and  $\L^{k_j}:C_j(0,1,1)\to C_{j+1}(0,1,1)$ admits, w.r.t. cone
metrics, a contraction coefficient less than:
$$
 \max\left(\frac15, \frac{v_{q(j)}}{v_{q(j+1)}} \right) =: \gamma_j.
$$
\end{prop}
Let us prove that if $   \L^{k_j} C_j(0,1,1) \subset    C_{j+1}\left(\frac45,\frac15,   \max\left(\frac15,\frac{v_{q(j)}}{v_{q(j+1)}} \right)   \right) $ then we have the announced estimation on the contraction rate.\\
This will follow if we prove that for all $f,g\in
C_{j+1}(4/5,1/5,\max(1/5,v_{q(j)}/v_{q(j+1)}))$ with  the
normalization $E_\muh(f)=E_\muh(g)=1$ we have:
$$
   \alpha f-g\in C_{j+1}(0,1,1)
$$
for
\begin{equation} \label{eqalpha}
\alpha = \max\left( \frac{1+D^{-1}}{1-D^{-1}},
\frac{1+v_{q(j)}/v_{q(j+1)}}{1-v_{q(j)}/v_{q(j+1)}} \right).
\end{equation}
\medbreak Indeed, in that case, we have that the diameter
$\Gamma_j$ of $\L^{k_j} C_j(0,1,1)$ into $C_{j+1} (0,1,1)$ is less
than
$$
2 \log \max\left( \frac{1+D^{-1}}{1-D^{-1}},
\frac{1+v_{q(j)}/v_{q(j+1)}}{1-v_{q(j)}/v_{q(j+1)} }\right).
$$
and then
$$
\tanh \frac{\Gamma_j}4 \leq \max \left( \frac1D
\/,\frac{v_{q(j)}}{v_{q(j+1)}} \right) = \max\left(\frac15,
\frac{v_{q(j)}}{v_{q(j+1)}} \right)  \/.
$$
The upper bound in the cone condition (1) for $\alpha f-g$ is, for
all $P\in\Q$,
$$
    \alpha \geq         \frac{6 E_\muh(g) - E_\muh(g|P)} {6 E_\muh(f) - E_\muh(f|P)}.
$$
The right hand side is bounded by:
$$
  \frac{6 - 0}{6-\frac65} = \frac54      \leq 6/4.
$$
The lower bound in this condition is, for all $P\in\Q$,
$$
 \alpha \geq \frac{ E_\muh(g|P) }{ E_\muh(f|P) }.
$$
The right hand side is upper bounded by:
$$
  \frac{6/5 \cdot E_\muh(g)}{4/5 \cdot E_\muh(f)} = \frac 64.
$$
Thus, both bounds in condition (1) are implied by eq.
(\ref{eqalpha}). \medbreak \noindent The cone condition (2) is,
for all $x,y\in\Delta$ with $\BP(x)=\BP(y)$
$$
  \alpha \geq \frac{ 12 E_\muh(g) + |g(x)-g(y)|}   { 12 E_\muh(f) - |f(x)-f(y)|}
$$
The right hand side is bounded by:
$$
  \frac{1+D^{-1}}{1-D^{-1}} = 6/4.
$$
Thus, condition (2) is implied by eq. (\ref{eqalpha}). \medbreak
\noindent The cone condition (3) is implied by eq. (\ref{eqalpha})
as can be seen by practically identical computations. \medbreak
\noindent The cone condition (4) is satisfied iff, for all $x\in
P_\infty\cap\Delta_\ell$, $\ell>q(j+1)$,
$$
 \alpha \geq \frac{ 90 v_{q(j+1)} E_\muh(g) + |g(x)|}{ 90 v_{q(j+1)} E_\muh(f) - |f(x)|}.
$$
But the right hand side is bounded by:
$$
  \frac{1+(v_{q(j)}/v_{q(j+1)})}{1-(v_{q(j)}/v_{q(j+1)})}.
$$
Thus, condition (4) is implied by eq. (\ref{eqalpha}) and this
concludes the proof that the claim implies the stated contraction
coefficient.
\subsection{Contraction of the first condition}
$f$ is an arbitrary function in $C_j(0,1,1)$ for the remainder of
section \ref{sec:contraction}. \medbreak \noindent Let $P\in\Q$.
We first prove the lower bound:
$$
\begin{aligned}
   E_\muh(\L^{k_j}f|P) &= \frac1{\muh(P)} \int_P \L^{k_j}f \, d\muh     = \frac1{\muh(P)} \int_\Delta 1_P \cdot \L^{k_j}f \, d\muh \\
      & = \frac1{\muh(P)} \int_\Delta 1_P\circ F^{k_j} \cdot f \,d\muh       = \frac1{\muh(P)} \int_{F^{-k_j}P} f \,d\muh \\
  &\geq \sum_{P'\in\Q_1} \frac1{\muh(P)} \int_{F^{-k_j}P\cap P'} f \, d\muh  + \frac1{\muh(P)} \int_{F^{-k_j}P\cap P_\infty} f \, d\muh \\
 & \geq \sum_{P'\in\Q_1} \frac{\muh(F^{-k_j}P\cap P')\muh(P')}{\muh(P)\muh(P')}  \biggl\{ E_\muh(f|P') -12(D+1)R_0(s) E_\muh(f) \biggr\} \\
 & \qquad   - \sum_{\ell\geq0} \frac{\muh(F^{-k_j}P\cap P_\infty\cap\Delta_\ell)}{\muh(P)}  \cdot 90 v_{\min(\ell,q(j))} E_\muh(f),
 \end{aligned}
$$
using $\diam_{d_j}(\Q_1)\leq (D+1) R_0(s)$, conditions (2)-(4). We
continue (obviously:  $v_\ell \geq v_{\min(\ell,q(j))}$):
$$
\begin{aligned}
 E_\muh(\L^{k_j}f|P)    & \geq \sum_{P'\in\Q_1} \frac{\muh(F^{-k_j}P\cap P')}{\muh(P)\muh(P')}    \muh(P')  \biggl\{ E_\muh(f|P') -12(D+1)R_0(s) E_\muh(f) \biggr\} \\
 & \qquad\qquad - \frac{\muh_v(F^{-k_j}P\cap P_\infty)}{\muh(P)\muh_v(P_\infty)} \muh_v(P_\infty)   \cdot 90  E_\muh(f)
\end{aligned}
$$
$$
\begin{aligned}
   &\geq \sum_{P'\in\Q_1} \frac78 \muh(P')E_\muh(f|P')   - \sum_{P'\in\Q_1} \frac78 \muh(P') \cdot 12 (D+1) R_0(s)   E_\muh(f)     \\
& \qquad\qquad-\frac98 \muh_v(P_\infty) \cdot 90 E_\muh(f) \\
  & \geq \frac78 \biggl\{ E_\muh(f)-\int_{P_\infty} f\, d\muh \biggr\}   -12\frac98 (D+1) R_0(s) E_\muh(f) \\
   & \qquad\qquad - 90 \frac98 \muh_v(P_\infty) E_\muh(f). \\
 \end{aligned}
$$
Observe that:
 $$
\begin{aligned}
 \int_{P_\infty} f \, d\muh &= \int_{P_\infty} \frac{f}{v} \, d\muh_v  \leq \sum_{\ell\geq0} \muh_v(P_\infty\cap\Delta_\ell) \cdot     90 \frac{ v_{\min(\ell,q(j))} }{ v_\ell } E_\muh(f) \\
  & \leq 90 \muh_v(P_\infty) E_\muh(f).
\end{aligned}
$$
Hence,
$$
\begin{aligned}
E_\muh(\L^{k_j}f|P)   & \geq \left\{ \frac78 -90 \left(\frac98+\frac78\right) \muh_v(P_\infty)    - 12 \frac98 (D+1) R_0(s) \right\} E_\muh(f)\\
   & \geq \frac45 E_\muh(f).
 \end{aligned}
$$
Similarly, we get the upper bound,
$$
\begin{aligned}
  E_\muh(\L^{k_j}f|P) & \leq \left\{ \frac98 +90 \left(\frac98+\frac78\right) \muh_v(P_\infty)     +12 \frac98 (D+1) R_0(s) \right\} E_\muh(f)\\
 & \leq 6D^{-1} E_\muh(f).
  \end{aligned}
$$
\subsection{Contraction of the second condition}\label{contrac2}
Let $x,y\in\Delta_\ell$ with $\BP(x)=\BP(y)$. First assume that $\ell\geq k_j$. Setting \\
$x^-=(x,\ell-k_j)$, $y^-=(y,\ell-k_j)\in\Delta$ (with a slight
abuse of notation), we have
$$
\begin{aligned}
|\L^{k_j}f(x) - \L^{k_j}f(y)| &= |f(x^-)-f(y^-)|  \leq 12 d_j(x^-,y^-) E_\muh(f) \\
&= 12 R_j(s(x,y)+k_j)  \leq 12 D^{-1} d_{j+1}(x,y).
\end{aligned}
$$
Now assume that $\ell<k_j$. We have $|\L^{k_j}f(x) - \L^{k_j}f(y)|
= |\L^rf(x^0)-\L^rf(y^0)|$ with $r=k_j-\ell$ and $x^0=(x,0)$,
$y^0=(y,0)$ (with the same abuse). Hence it is enough to bound
$|\L^rf(x)-\L^rf(y)|$ for $r\leq k_j$ and $x,y\in\Delta_0$ with
$\BP(x)=\BP(y)$.   As $\BP(x)=\BP(y)$, the pre-images by $F^r$ of
$x$ and $y$ can be paired (i.e., to each pre-image $x'$ of $x$
corresponds a pre-image $y'$ of $y$ defined by the same inverse
branch). Thus,
\begin{eqnarray*}
  |\L^rf(x)-\L^rf(y)|  &\leq \sum_{x'\in F^{-r}x} \left| \frac{f(x')}{JF^r(x')} -  \frac{f(y')}{JF^r(y')}\right|\\
   &\leq \sum_{x'\in F^{-r}x} \frac1{JF^r(x')} | f(x') - f(y') | + \\
    & \qquad\qquad  + \sum_{x'\in F^{-r}x} |f(y')| \frac1{JF^r(x')}   \left| \frac{JF^r(x')}{JF^r(y')} - 1 \right| \\
  &\leq 12  R_j(r+s(x,y)) E_\muh(f) \\
   &  + C d_0(x,y) \left( \sum_{x'\in P_\infty \atop x'\in F^{-r}x}     \frac{|f(y')|}{JF^r(x')}       + \sum_{x'\notin P_\infty \atop x'\in F^{-r}x}         \frac{|f(y')|}{JF^r(x')}  \right)
\end{eqnarray*}
recall $\L1=1$ and $C$ is defined in (A.III). We have
$$
\begin{aligned}
  \sum_{x'\in P_\infty \atop x'\in F^{-r}x}  & \leq \sum_{\ell\geq0} \sum_{x'\in P_\infty\cap\Delta_\ell \atop x'\in F^{-r}x}            90 v_{\min(\ell,q(j))} \frac{1}{JF^r(x')} E_\muh(f)\\
   &  \leq 180 \frac{K}{\eta} \muh_v(P_\infty) E_\muh(f)
\end{aligned}
$$
where $K$ is given by the bounded distortion and $\eta$ by the large image property (we have used that $\int_{F^rC_r(x')} (JF^r)^{-1} \,d\muh  =\muh(C_r(x'))\geq\eta/K (JF^r(x'))^{-1}$).  \\
We also have
$$
\begin{aligned}
  \sum_{x'\notin P_\infty \atop x'\in F^{-r}x}   & \leq  \biggl\{ E_\muh(f|\Q)(x')+12R_j(s)E_\muh(f) \biggr\} \\
    & \leq  (12R_j(s)+1)E_\muh(f).
  \end{aligned}
$$
Hence,
$$
\begin{aligned}
 |\L^{k_j}f(x)-\L^{k_j}f(y)|  & \leq \left\{  12R_j(k_j-\ell+s(x^0,y^0))+d_0(x^0,y^0)MC(12R_j(s) \right.\\
& \left. +1   + 180 \frac{K}{\eta}\muh_v(P_\infty) ) \right\} E_\muh(f)\\
   & \leq 12\left\{ R_j(k_j+s(x,y))+R_0(s(x,y))(MCR_j(s)+1/12    \right.\\
   &\left.   + 15 \frac{K}{\eta}\muh_v(P_\infty)) \right\} E_\muh(f)
\end{aligned}
$$
Now, $CR_j(s)+ 15 (K/\eta) \muh_v(P_\infty) <1/2$ so that
 $$
|\L^{k_j}f(x)-\L^{k_j}f(y)| \leq 12D^{-1} R_{j+1}(s(x,y))
E_\muh(f)      = 12D^{-1} d_{j+1}(x,y) E_\muh(f).
$$
\subsection{Contraction of the third condition}
Let $x\in\Delta_\ell\cap P_\infty$. First assume $0\leq \ell\leq
k_j$. As for the second condition it is enough to consider
$\L^rf(x)$ with $0\leq r\leq k_j$ and $x\in\Delta_0$. We have
$$
\begin{aligned}
   |\L^rf(x)| &\leq \sum_{x'\in F^{-r}x} \frac{1}{JF^r(x')}|f(x')| \\
   &\leq \sum_{x'\in F^{-r}x \atop x'\notin P_\infty}  \frac{1}{JF^r(x')}\left( 12R_j(s)E_\muh(f)+6E_\muh(f) \right) \\
     &\qquad +  \sum_{\ell\geq0}\sum_{x'\in F^{-r}x \atop x'\in P_\infty\cap\Delta_\ell}   \frac{1}{JF^r(x')} 90 v_{\min(\ell,q(j))} E_\muh(f) \\
     & \leq  (2R_j(s)+6) E_\muh(f)    +  \sum_{\ell\geq0}\sum_{x'\in F^{-r}x \atop x'\in P_\infty\cap\Delta_\ell}        \frac{K}{\eta} \muh(C_r(x')) \cdot 90 v_{\min(\ell,q(j))} E_\muh(f) \\
  & \leq  (12R_j(s)+6) E_\muh(f)      +  \sum_{\ell\geq0}\sum_{x'\in F^{-r}x \atop x'\in P_\infty\cap\Delta_\ell}       \frac{K}{\eta} \muh_v(C_r(x')) \cdot 90 v_{\min(\ell,q(j))} E_\muh(f) \\
    & \leq \left( (12R_j(s)+6)+90\frac{K}{\eta}\muh_v(P_\infty) \right) E_\muh(f) \\
& \leq 90 D^{-1} E_\muh(f)
\end{aligned}
$$
Now assume $k_j\leq\ell\leq q(j+1)=q(j)+k_j$ and let
$x^-=(x,\ell-k_j)$. We have:
$$
\begin{aligned}
    |\L^{k_j}f(x)| &= |f(x^-)| \leq 90 v_{\ell-k_j} E_\muh(f) \\
    & \leq 90 \frac{v_{\ell-k_j}}{v_\ell} v_\ell E_\muh(f) \\
    & \leq 90 \frac{v_{q(j)}}{v_{q(j+1)}} v_\ell E_\muh(f).
 \end{aligned}
$$
using that $\ell\mapsto v_\ell/v_{\ell+k}$ is increasing for any
$k$.
\\
We  get the claimed contraction by
$\max(1/5,\frac{v_{q(j)}}{v_{q(j+1)}})$.
\subsection{Contraction of the fourth condition}
Finally we take $x\in P_\infty\cap\Delta_\ell$ with $\ell>q(j+1)$.
We have
$$
\begin{aligned}
   |\L^{k_j}f(x)| &= |f(x^-)| \leq 90 v_{q(j)} E_\muh(f) \\
     & \leq 90 \frac{v_{q(j)}}{v_{q(j+1)}} v_{q(j+1)} E_\muh(f).
\end{aligned}
$$
and this gives the contraction by $\frac{v_{q(j)}}{v_{q(j+1)}}$.
\section{Conclusion}\label{sec:conclusion}
To conclude the proof of Theorem \ref{theotour}, we need to derive
from the projective metric bound obtained above, bounds on the
correlations. The following lemma is standard when using
Birkhoff's cones (see \cite{KMS} page 687, \cite{M} Lemmas
3.9-3.10). Let $||| \ |||_j$ be the norm on bounded functions
defined by:
\begin{eqnarray*}
||| f |||_j& =& \max \left[ \max(90 v_{q(j)} \/, 12 D R_0(s)+6) \left\vert \int\limits_\Delta f d\muh \right\vert \/,\right. \\
& &\left. \ \sup_{ P \in \Q} \muh(P)^{-1} \left| \int\limits_P f
d\muh \right| \/, \ \Vert f \Vert_\infty \right] \/.
\end{eqnarray*}
\begin{lem}\label{pptes}
The  norms $||| \ |||_j$  and the homogeneous form $f \mapsto \muh (f)$ are adapted to the cones $C_j(0\/,1\/,1)$\\
For any $f \in {\bf L}(d_0)$, there exists $R(f) >0$ such that $f
+ R(f) \un \in C_0 (0,1,1)$  and $R(f) \leq C \Vert f
\Vert_{\LL}$.
\end{lem}
\begin{sketchdem}
It is clear that the  homogeneous form $f \mapsto \muh (f)$ is
adapted. To prove that $||| \ |||_j$ is also adapted, let us
consider $f$ and $g$ such that $f+g$ and $f-g$ are in
$C_j(0\/,1\/,1)$. The first condition in the definition of the
cone gives:
$$
\forall \ P \in \Q \/, \left| \frac1{\muh(P)} \int\limits_P g
d\muh \right| \leq \frac1{\muh(P)} \int\limits_P f d\muh  \/,
$$
and
$$
\left|\int_\Delta g d\muh  \right|\leq \int_\Delta f d\muh \/.
$$
The last three conditions give
$$
\Vert g \Vert_\infty \leq \max [90 v_{q(j)} \/, 12 D R_0(s)+6] \/
\int_\Delta f d\muh.
$$
Hence, we have $||| g |||_j \leq |||f|||_j$.\\
To prove the second point of the lemma, we may assume that
$f\geq0$. To have that  $f + R(f)\un \in C_0 (0,1,1)$, it suffices
that:
\begin{itemize}
\item $\forall P \in \Q$, $R(f) \geq
\frac{\frac1{\muh(P)}\int_\Delta P f d\muh  -6\int_\Delta f
d\muh}{5}$, so that condition 1. is satisfied, \item $R(f) \geq
\frac{L(f)}{12}$,  so that condition 2. is satisfied, \item $R(f)
\geq \frac{\sup f}{90 v_{q(0)} -1}$,  so that conditions 3 and 4
are  satisfied,
\end{itemize}
so we may choose $R(f) \leq \cte ||f||_\LL$.
\end{sketchdem}
\ \\
\noindent Let us conclude the proof of Theorem \ref{theotour}.\\
Let  $f \in C_0(0,1,1)$, By Proposition \ref{image}, for any
$\ell$, $\L^{k_1+ \cdots + k_\ell } f$ and $\L^{k_1+ \cdots +
k_\ell } \un = \un$ belong to $C_\ell$ (we remark that $\un \in
C_0(0,1,1)$). Applying $\ell -1$ times Theorem \ref{projm1}, we
get:
$$
\delta_{C_{\ell }} (\L^{k_1+ \cdots + k_\ell } f , \un) \leq
\prod_{j=2}^{\ell} \gamma_j \cdot \delta_1 (\L^{k_1}f ,\un) \leq
\prod_{j=2}^{\ell} \gamma_j \cdot \Gamma_{1}\/,
$$
where $\gamma_j=$ is given by Proposition \ref{image}. Since the
norm $|||  \ |||_j$ is adapted to the cones $C_j (0,1,1)$ and is
greater than  the uniform  norm,  Theorem \ref{projm2} gives for
$f$  in $C_0(a,b,c)$ with $\muh (f) =1$,
$$
 \Vert \L^{k_1 + \cdots k_\ell} f - \un \Vert_\infty \leq \cte \/\prod_{j=2}^{\ell} \gamma_j  \/.
$$
For $n \in \NN$, let $\ell(n)$ be defined by:
$$
n =k_1 + \cdots k_{\ell(n)} +r \ \mbox{with} \ r < k_{\ell (n) +1}
\/,
$$
we have
$$
 \Vert \L^{n} f -  \un \Vert_\infty \leq \Vert \L^r \un \Vert_\infty \ \Vert \L^{k_1 + \cdots + k_\ell} f - \un \Vert_\infty \leq \cte \/\left(\prod_{j=2}^{\ell(n)} \gamma_j \right) \/.
$$
For any function $f \in \LL$, applying the above inequality to
$\frac{f+R(f)\un}{\muh(f) + R(f)}$ gives
$$
 \Vert \L^n f - \muh (f) \Vert_\infty \leq \cte \cdot \left(\prod_{j=2}^{\ell(n)} \gamma_j\right)  \Vert f \Vert_\LL \/.
$$
The decay of correlations follows: for $f \in L$ and $g \in
L^{1}(\muh)$,
\begin{eqnarray*}
\left\vert \int\limits_\Delta g \circ F^n \cdot f d\muh  - \int\limits_\Delta f d\muh \int\limits_\Delta g d\muh \right\vert &=&\left\vert \int\limits_\Delta g [ \L^n (f) -  \muh(f) ] d\muh \right\vert \\
&\le& \cte \cdot \left( \prod_{j=2}^{\ell(n)} \gamma_j \right)
\Vert f \Vert_\LL \/ \Vert g \Vert_1 \/.
\end{eqnarray*}
Set $\prod_{j=2}^{\ell(n)} \gamma_j =u_n $.
\\
\\
We now prove that $u_n$ has the announced behavior for
exponential, stretched exponential or polynomial sequences
$\ItOsc_n$  and $\nuh(\Delta_n)$.
\subsubsection*{Estimating $u_n$}
To estimate the rate of mixing $u_n$, we have to analyze the
asymptotic behavior of the sequences $k_j$ and $\gamma_j$. Recall
that:
$$
 k_{\ell+1}=\inf\{k\geq k_0: R_\ell(s+k)\leq \frac{R_0(s)}{ D}  \}
$$
and observe that an easy induction  gives:
 \begin{equation} \label{eq:Rell}
R_\ell (n) = D ^\ell R_0 (k_1 + \cdots + k_\ell +n) +
\sum_{i=1}^{\ell}  D^{\ell +1 -i} R_0(k_{i+1} + \cdots + k_\ell
+n) \/.
 \end{equation}
Recall $R_0(m)=\sum_{k>m} \omega_k$.\\
\\
From these remarks we get the following results.
\begin{itemize}
\item  If there exists $0<\rho<1$ such that $\om_n = O (\rho^n) $
then $R_0(n)= O(\rho^n)$ and  one can take  $k_j=p$ provided $p$
is such that $\rho^p< \frac1{D(D+1)}$, i.e. $p
>\displaystyle\frac{-\log[D(D+1)]}{\rho}$. Assume also that $\nuh
(\Delta_n) =O(\alpha^n)$  for some $0< \alpha < 1$. Then one may
choose $v_n={\alpha'}^{-n}$ provided $0<\alpha<{\alpha'}<1$. We
have: $\gamma_j  = \max\left(\frac1D \/, {\alpha'}^p
\right)=:\kappa < 1$ and $u_n = \kappa^{\frac{n}p}$.
\item If there exists $\alpha >1$ such that $\om_n = O (n^{-\alpha})$ then  $R_0(n) = O(n^{-\alpha+1})$  and $R_\ell (s+k_{\ell+1}) =  O\left(R_0 (s+ k_{\ell +1}) \cdot D^\ell\right) = O\left(\frac{D^\ell}{k_{\ell+1}^{\alpha-1}}\right)$ so,  if $k_{\ell+1} \sim \cte \/ D^{\frac{\ell+1}{\alpha -1}}$, it satisfies $R_\ell(s+k_{\ell+1}) \leq \frac{1}{DR_0(s)}$.  So $\ell (n) \sim (\alpha -1)\frac{\log n}{\log  D} + \cte$. Assume also that $\nuh (\Delta_n) =O(n^{-\beta})$ for some $\beta> 1$. Then one may choose $v_n= n^\gamma$ provided $0< \gamma< \beta-1$. We have: $\gamma_j=\max\left(\frac1D \/, D^{-\frac{\gamma}{\alpha-1}} \right)$.\\
If $\beta \leq \alpha$ then we can choose $\gamma <\beta-1 \leq\alpha-1$ and then $\gamma_j = D^{-\frac{\gamma}{\alpha-1}}$ and $u_n=O(n^{-\gamma})$.\\
 If $\alpha < \beta$ then we can choose $\beta-1>\gamma >\alpha -1$ and then $\gamma_j=\frac1D$ and $u_n = O(n^{-\alpha +1})$.  \\
Finally, we get that $u_n =O(n^{-\min[\alpha-1 \/,
\beta-1-\eps]})$, for all $\eps >0$.
\item If $\om_n=O(e^{-n^\alpha})$ and
$\nuh(\Delta_n)=O(e^{-n^\beta})$ for some $0<\alpha,\beta<1$, set
$k_\ell:= \ell^{\frac1\alpha-1}$. We obtain
$q(\ell)=k_1+\dots+k_\ell \sim \ell^{\frac1\alpha}$ so
$\ell(n)\sim n^{\alpha}$. An easy estimate gives that $R_0(m) \leq
e^{-m^{\alpha-\eps}}$ for $\eps>0$ and all large $m$. Now, by eq.
(\ref{eq:Rell}),
$$
  R_\ell(s+k_{\ell+1}) \leq (\ell+1) D^\ell R_0(k_{\ell+1}) \leq (\ell+1) D^\ell e^{-\ell^{(\alpha-\eps)(\frac1\alpha-1)}} \to 0 \  \text{as }\ell\to\infty.
$$
Hence, $R_\ell(s+k_{\ell+1})\leq R_0(s)/D$ is satisfied and the
choice of $k_\ell$ is correct for all large $\ell$. Let us
compute:
\begin{eqnarray*}
\gamma_\ell& =& \max\left\{ \frac1D, \frac{q(\ell+1)^2}{q(\ell)^2} \exp\left(\ell^{\frac\beta\alpha}- (\ell+1)^{\frac\beta\alpha} \right) \right\}  \\
&  =& \max\left\{ \frac1D,
\left(1+\frac1\ell\right)^{\frac2\alpha} \exp
\left(-\ell^{\frac\beta\alpha}\cdot \frac\beta\alpha \frac1\ell +
\dots \right) \right\}.
\end{eqnarray*}
If $\beta>\alpha$, then the second term of the above maximum goes
to zero and therefore $\gamma_\ell=\frac1D$ for large $\ell$. We
compute the contraction coefficient at time $n$:
$$
   u_n = D^{-\ell(n)} = D^{- C n^\alpha} = e^{- C' n^\alpha}\leq e^{- n^{\alpha-\eps}}.
$$
If $\beta<\alpha$, then the second term of the above maximum goes
to one and therefore sets the value of $\gamma_\ell$ for large
$\ell$. We compute:
$$
   u_n = \prod_{j=1}^{\ell(n)} e^{-j^\frac\beta\alpha} \leq e^{-\ell(n)^\frac\beta\alpha} = e^{- C n^\beta}   \leq e^{-n^{\beta-\eps}}.
$$
If $\alpha=\beta$, then we change, for instance, $\alpha$ to
$\alpha'>\beta$, arbitrarily close, and apply the previous case.
\end{itemize}

{\small \noindent {\em e-mail:} {\tt buzzi@math.polytechnique.fr  }\\
\noindent Centre de Math{\'e}matiques de l'Ecole Polytechnique\\
\noindent U.M.R. 7640 du C.N.R.S.\\
\noindent Ecole Polytechnique\\
\noindent 91128 Palaiseau Cedex / FRANCE  \ \\
\noindent {\em e-mail:} {\tt vmaume@topolog.u-bourgogne.fr}\\
\noindent {\sc Laboratoire de Topologie }\\
\noindent Universit{\'e} de Bourgogne\\
\noindent    9,Avenue Alain Savary -   B.P. 47870\\
\noindent    21078 Dijon Cedex FRANCE}

\end{document}